%
\input amstex
\input amsppt.sty
\magnification=\magstep1
\hsize = 6.5 truein
\vsize = 8 truein
\NoBlackBoxes
\TagsAsMath
\NoRunningHeads
%
%
\define\op#1{\operatorname{\fam=0\tenrm{#1}}} 
\redefine\qed{{\unskip\nobreak\hfil\penalty50\hskip2em\vadjust{}\nobreak\hfil
$\square$\parfillskip=0pt\finalhyphendemerits=0\par}}

	\let	\< = \langle
	\let	\> = \rangle
	\define		 \a		{\alpha}
	\define		 \g		{\frak g}
\topmatter
\title Metrics on States from Actions of Compact Groups 
\endtitle
\author Marc A. Rieffel
\endauthor
\address\hskip-\parindent
Department of Mathematics \newline University of California
\newline Berkeley, California 94720-3840
\endaddress
\email rieffel\@math.berkeley.edu
\endemail
\date	July 14, 1998
\enddate
\thanks The research reported here was supported in part by National Science
Foundation\newline Grant DMS--96--13833.
\endthanks
\abstract{Let a compact Lie group act ergodically on a unital $C^*$-algebra
$A$. We consider several ways of using this structure to define metrics on the
state space of $A$. These ways involve length functions, norms on the Lie
algebra, and Dirac operators. The main thrust is to verify that the
corresponding metric topologies on the state space agree with the weak-$*$
topology. }
\endabstract
\endtopmatter
\document
\baselineskip=16 pt

Connes \cite{C1, C2, C3} has shown us that Riemannian metrics on
non-commutative spaces ($C^*$-algebras) can be specified by generalized
Dirac operators. Although in this setting there is no underlying manifold
on which one then obtains an ordinary metric, Connes has shown that one
does obtain in a simple way an ordinary metric on the state space of the
$C^*$-algebra, generalizing the Monge-Kantorovich metric on probability
measures \cite{Ra} (called the ``Hutchinson metric'' in the theory of
fractals \cite{Ba}).

But an aspect of this matter which has not received much attention so far
\cite{P} is the question of when the metric topology (that is, the topology
from the metric coming from a Dirac operator) agrees with the underlying
weak-$*$ topology on the state space. Note that for locally
compact spaces their topology agrees with the weak-$*$ topology coming from
viewing points as linear functionals (by evaluation) on the algebra of
continuous functions vanishing at infinity.

In this paper we will consider metrics arising from actions of compact
groups on $C^*$-algebras. For simplicity of exposition we will only deal
with ``compact" non-commutative spaces, that is, we will always assume that
our $C^*$-algebras have an identity element. We will explain later what we
mean by Dirac operators in this setting (section 4). In terms of this, a
brief version of our main theorem is:

\proclaim{Theorem 4.2} Let $\a$ be an ergodic action of a compact Lie group
$G$ on a unital $C^*$-algebra $A$, and let $D$ be a corresponding Dirac
operator. Then the metric topology on the state space of $A$ defined by the
metric from $D$ agrees with the weak-$*$ topology.
\endproclaim

An important case to which this theorem applies consists of the
non-commutative tori \cite{Rf}, since they carry ergodic actions of
ordinary tori \cite{OPT}. The metric geometry of non-commutative tori has
recently become of interest in connection with string theory \cite{CDS, RS,
S}.

We begin by showing in the first section of this paper that the mechanism
for defining a metric on states can be formulated in a very rudimentary
Banach space setting (with no algebras, groups, or Dirac operators). In
this setting the discussion of agreement between the metric topology and the
weak-$*$ topology takes a particularly simple form.

Then in the second section we will see how length functions on a compact
group directly give (without Dirac operators) metrics on the state spaces
of $C^*$-algebras on which the group acts ergodically. We then prove the
analogue in this setting of the main theorem stated above.

In the third section we consider compact Lie groups, and show how
norms on the Lie algebra directly give metrics on the state space.
We again prove the corresponding analogue of our main theorem.

Finally, in section 4 we use the results of the previous sections to prove
our main theorem, stated above, for the metrics which come from Dirac
operators.

It is natural to ask about actions of non-compact groups. Examination of
\cite{Wv4} suggests that there may be very interesting phenomena there. The
considerations of the present paper also make one wonder whether there is
an appropriate analogue of length functions for compact quantum groups
which might determine a metric on the state spaces of $C^*$-~algebras on
which a quantum group acts ergodically \cite{Bo, Wn}. This would be
especially interesting since for non-commutative compact groups there is
only a sparse collection of known examples of ergodic actions \cite{Ws},
whereas in \cite{Wn} a rich collection of ergodic actions of compact
quantum groups is constructed. Closely related is the setting of ergodic
coactions of discrete groups \cite{N, Q}. But I have not explored any of
these possibilities.

I developed a substantial part of the material discussed in the present
paper during a visit of several weeks in the Spring of 1995 at the Fields
Institute. I am appreciative of the hospitality of the Fields Institute,
and of George Elliott's leadership there. But it took trying to present
this material in a course which I was teaching this Spring, as well as
benefit from
\cite{P, Wv1, Wv2, Wv3, Wv4}, for 
me to find the simple development given here.

\bigskip\bigskip\noindent{\bf 1. Metrics on states}

\medskip
\noindent
Let $A$ be a unital $C^*$-algebra. Connes has shown \cite{C1, C2, C3}
that an appropriate way to specify a Riemannian metric in this
non-commutative situation is by means of a spectral triple. This consists
of a representation of $A$ on a Hilbert space $\Cal H$, together with an
unbounded self-adjoint operator $D$ on $\Cal H$ (the generalized Dirac
operator), satisfying certain conditions. The set $\Cal L(A)$ of Lipschitz
elements of $A$ consists of those $a\in A$ such that the commutator $[D,a]$
is a bounded operator. It is required that $\Cal L(A)$ be dense in $A$. The
Lipschitz semi-norm, $L$, is defined on $\Cal L(A)$ just by the operator
norm $L(a)=\|[D,a]\|$.

Given states $\mu$ and $\nu$ of $A$, Connes defines the distance between
them, $\rho(\mu,\nu)$, by
$$
\rho(\mu,\nu)=\sup\{ |\mu(a)-\nu(a)|: a\in\Cal L(A), \ L(a)\Cal \ \leq \
1\} \ .
\tag 1.1
$$
(In the absence of further hypotheses it can easily happen that
$\rho(\mu,\nu)=+\infty$. For one interesting situation where this
sometimes happens see the end of the discussion of the second
example following axiom 4' of \cite{C3}.)

The semi-norm $L$ is an example of a general Lipschitz semi-norm, that is
\cite{BC, Cu, P, Wv1, Wv2}, a semi-norm
$L$ on a dense subalgebra $\Cal L$ of $A$ satisfying the Leibniz property:
$$
L(ab) \ \Cal \leq \ L(a)\|b\| +\|a\| L(b) \ . \tag 1.2 $$
Lipschitz norms carry some information about differentiable structure
\cite{BC, Cu}, but not nearly as much as do spectral triples. But it is
clear that just in terms of a given Lipschitz norm one can still define a
metric on states by formula (1.1).

However, for the purpose of understanding the relationship between the
metric topology and the weak-$*$ topology, we do not need the Leibniz
property (1.2), nor even that $A$ be an algebra. The natural setting for
these considerations seems to be the following very rudimentary one. The
data is: $$
{\text{A normed space $A$, with norm $\| \ \|$, over either $\Bbb C$ or
$\Bbb R$.}} \tag 1.3a
$$
$$
{\text{A subspace $\Cal L$ of $A$, not necessarily closed.}} \tag 1.3b $$
$$
{\text{A semi-norm $L$ on $\Cal L$.}} \tag 1.3c $$
$$
\aligned
&\text{A continuous (for } \| \ \|) \text{ linear functional, } \eta,
\text{ on } {\Cal K}= \{a\in\Cal L: L(a)=0\} \\
&\text{with }\|\eta\|=1 \text{. (Thus, in particular, we require} \ \ 
{\Cal K}\neq \{0\}\ .)
\endaligned \tag 1.3d
$$

Let $A'$ denote the Banach-space dual of $A$, and set $$
S=\{\mu\in A':\mu=\eta \ {\text{on}} \ \Cal K, \ {\text{and}} \ \|\mu\|=1\} \ .
$$
Thus $S$ is a norm-closed, bounded, convex subset of $A'$, and so is
weak-$*$ compact. In general $S$ can be quite small; when $A$ is a Hilbert
space $S$ will contain only one element. But in the applications we have in
mind $A$ will be a unital $C^*$-algebra, $\Cal K$ will be the
one-dimensional subspace spanned by the identity element, and $\eta$ will
be the functional on $\Cal K$ taking value 1 on the identity element. Thus
$S$ will be the full state-space of $A$. (That $\Cal K$ will consist only
of the scalar multiples of the identity element in our examples will follow
from our
ergodicity hypothesis. We treat the case of general $\Cal K$ here because
this clarifies slightly some issues, and it might possibly be of eventual
use, for example in non-ergodic situations.)

We do not assume that $\Cal L$ is dense in $A$. But to avoid trivialities
we do make one more assumption about our set-up, namely:
$$
{\Cal L} \ {\text{ \ separates the points of \ }} \ S. \tag 1.3e $$
This means that given $\mu,\nu\in S$ there is an $a\in\Cal L$ such that
$\mu(a)\neq\nu(a)$.
(Note that for $\mu\in S$ there exists $a\in\Cal L$ with $\mu(a)\neq 0$,
since we can just take an $a\in\Cal K$ such that $\eta(a)\neq 0$.)

With notation as above, let $\tilde\Cal L=\Cal L/\Cal K$. Then $L$ drops to
an actual norm on $\tilde\Cal L$, which we denote by $\tilde L$. But on
$\tilde\Cal L$ we also have the quotient norm from $\| \ \|$ on $\Cal L$,
which we denote by $\| \ \|^{\sptilde}$. The image in $\tilde\Cal L$ of
$a\in\Cal L$ will be denoted by $\tilde a$.

We remark that when $\Cal L$ is a unital algebra (perhaps dense in a
$C^*$-algebra), and when $\Cal K$ is the span of the identity element, then
the space of universal 1-forms $\Omega^1$ over $\Cal L$ is commonly
identified \cite{BC, Br, C2, Cu} with $\Cal L\otimes\tilde\Cal L$, and the
differential $d:\Cal L\to\Omega^1$ is given by $da=1\otimes\tilde a$. Thus
in this setting our $\tilde L$ is a norm on the space of universal
1-coboundaries of $\Cal L$. The definition of $L$ which we will use in the
examples of section 3 is also closely related to this view.

On $S$ we can still define a metric, $\rho$, by formula (1.1), with $\Cal
L(A)$ replaced by $\Cal L$. The symmetry of $\rho$ is evident, and the
triangle inequality is easily verified. Since we assume that $\Cal L$
separates the points of $S$, so will $\rho$. But $\rho$ can still take the
value $+\infty$. We will refer to the topology on $S$ defined by $\rho$ as
the ``$\rho$-topology'', or the ``metric topology'' when $\rho$ is
understood.

It will often be convenient to consider elements of $A$ as (weak-$*$
continuous) functions on $S$. At times this will be done tacitly, but when
it is useful to do this explicitly we will write $\hat a$ for the
corresponding function, so that $\hat a(\mu)=\mu(a)$ for $\mu\in S$.

Without further hypotheses we have the following fact. It is closely
related to proposition 3.1a of \cite{P}, where metrics are defined in terms
of linear operators from an algebra into a Banach space.

\proclaim{1.4 Proposition} The $\rho$-topology on $S$ is finer than the
weak-$*$ topology.
\endproclaim

\demo{Proof}
Let $\{\mu_k\}$ be a sequence in $S$ which converges to $\mu \in S$ for the
metric $\rho$. Then it is clear from the definition of $\rho$ that
$\{\mu_k(a)\}$ converges to $\mu(a)$ for any $a \in {\Cal L}$ with $L(a)
\leq 1$, and hence for all $a \in {\Cal L}$.

This says that ${\hat a}(\mu_k)$ converges to ${\hat a}(\mu)$ for all $a
\in {\Cal L}$. But $\hat {\Cal L}$ is a linear space of weak@-$*$
continuous functions on $S$ which separates the points of $S$ by assumption
(and which contains the constant functions, since they come from any $a \in
{\Cal K}$ on which $\eta$ is not $0$). A simple compactness argument shows
then that $\hat {\Cal L}$ determines the weak-$*$ topology of $S$. Thus
$\{\mu_k\}$ converges to $\mu$ in the weak-$*$ topology, as desired.	\qed
\enddemo

\bigskip
There will be some situations in which we want to obtain information about
$(\Cal L,L)$ from information about $S$. It is clear that to do this $S$
must ``see" all of $\Cal L$. The convenient formulation of this for our
purposes is as follows. Let $\| \ \|_{\infty}$ denote the supremum norm on
functions on $S$. Let it also denote the corresponding semi-norm on $\Cal
L$ defined by $\|a\|_{\infty}=\|\hat a\|_{\infty}$. Clearly $\|\hat
a\|_{\infty}\leq \|a\|$ for $a\in\Cal L$.

\bigskip
\definition{1.5 \ Condition} The semi-norm $\| \ \|_{\infty}$ on $\Cal L$
is a norm, and it is equivalent to the norm $\| \ \|$, so that there is a
constant $k$ with
$$
\|a\|\leq k\|\hat a\|_{\infty} \qquad {\text{for}} \ \ a\in\Cal L. $$
\enddefinition
\bigskip

This condition clearly holds when $A$ is a $C^*$-algebra, $\Cal L$ is dense
in $A$, and $S$ is the state space of $A$, so that we are dealing with the
usual Kadison functional representation \cite{KR}. But we remark that even
in this case the constant $k$ above cannot always be taken to be 1 (bottom
of page 263 of \cite{KR}). This suggests that in using formula (1.1) one
might want to restrict to using just the self-adjoint elements of $\Cal L$,
since there the function representation is isometric. But more experience
with examples is needed.

We return to the general case. If we are to have the $\rho$-topology on $S$
agree with the weak-$*$ topology, then $S$ must at least have finite
$\rho$-diameter, that is, $\rho$ must be bounded. The following proposition
is closely related to theorem 6.2 of \cite{P}.

\proclaim{1.6 \ Proposition} Suppose there is a constant, $r$, such that
$$
\| \ \|^{\sptilde}\leq r\tilde L \ . \tag 1.7 $$
Then $\rho$ is bounded (by $2r$).

Conversely, suppose that Condition 1.5 holds. If $\rho$ is bounded, (say by
$d$), then there is a constant $r$ such that (1.7) holds (namely $r=kd$
where $k$ is as in 1.5).
\endproclaim

\bigskip
\demo{Proof} Suppose that (1.7) holds. If $a\in\Cal L$ and $L(a)\leq 1$,
then $\tilde L(\tilde a)\leq 1$ and so $\|\tilde a\|^{\sptilde}\leq r$.
This means that, given $\varepsilon > 0$, there is a $b\in\Cal K$ such that
$\|a-b\|\leq r+\varepsilon$. Then for any $\mu,\nu\in S$, we have, because
$\mu$ and $\nu$ agree on ${\Cal K}$, $$
|\mu(a)-\nu(a)| = |\mu(a-b)-\nu(a-b)|\leq \|\mu-\nu\| \ \|a-b\|\leq
2(r+\varepsilon) \ . $$
Since $\varepsilon$ is arbitrarily small, it follows that $|\mu(a)-\nu(a)|
\leq 2r$. Consequently $\rho(\mu,\nu)\leq 2r$.

Assume conversely that $\rho$ is bounded by $d$. Fix $\nu\in S$, and choose
$b\in\Cal K$ such that $\eta(b)=1$. Then for any $\mu\in S$ and any
$a\in\Cal L$ with $L(a)\leq 1$ we have $$
d\geq\rho(\mu,\nu)\geq |\mu(a)-\nu(a)| = |\mu(a-\nu(a)b)| \ . $$
Suppose now that Condition 1.5 holds. We apply it to $a-\nu(a)b$. Thus,
since $S$ is compact, we can find $\mu$ such that $$
\|a-\nu(a)b\|\leq k|\mu(a-\nu(a)b)| \ .
$$
Consequently $\|a-\nu(a)b\|\leq kd$, so that $\|\tilde a\|^{\sptilde}\leq
kd$. All this was under the assumption that $L(a)\leq 1$. It follows that
for general $a\in\Cal L$ we have $\|\tilde a\|^{\sptilde}\leq kd\tilde
L(\tilde a)$, as desired. \qed \enddemo

\bigskip
We now turn to the question of when the $\rho$-topology and the weak-$*$
topology on $S$ agree. The following theorem is closely related to theorem
6.3 of \cite{P}.

\proclaim{1.8 \ Theorem} Let the data be as in {\rm{(1.3a--e)}}, and let
$\Cal L_1=\{a\in\Cal L: L(a)\leq 1\}$. If the image of $\Cal L_1$ in $\Cal
L^{\sptilde}$ is totally bounded for $\| \ \|^{\sptilde}$, then the
$\rho$-topology on $S$ agrees with the weak-$*$ topology.

Conversely, if Condition 1.5 holds and if the $\rho$-topology on $S$ agrees
with the weak-$*$ topology, then the image of $\Cal L_1$ in $\Cal
L^{\sptilde}$ is totally bounded for $\| \ \|^{\sptilde}$.
\endproclaim

\demo{Proof} We begin with the converse,
so that we see why the total-boundedness assumption is natural.
If the $\rho$-topology gives the weak-$*$ t
opology on $S$, then $\rho$ must be bounded since $S$ is compact. Thus by
Proposition 1.6 there is a constant, $r_o$, such that $\| \ \|^{\sptilde}
\leq r_oL^{\sptilde}$, since we assume here that Condition 1.5 holds.
Choose $r > r_o$. Then
$\|a\|^{\sptilde} < r$ if $a\in\Cal L_1$. Consequently, if
we let
$$
\Cal B_r=\{a\in \Cal L: L(a)\leq 1 \ {\text{and}} \ \|a\|\leq r\} \ ,
$$
then the image of $\Cal B_r$ in $\Cal L^{\sptilde}$ is the
same as the image of
$\Cal L_1$. Thus it suffices to show that $\Cal B_r$ is totally bounded.

Let $a\in\Cal B_r$ and let $\mu,\nu\in S$. Then
$$
|\hat a(\mu)-\hat a(\nu)| = |\mu(a)-\nu(a)| \leq \rho(\mu,\nu) \ .
$$
Thus $(\Cal B_r)\ \hat{}$  can be viewed as a bounded family of functions on
$S$ which is equi-continuous for the weak-$*$ topology, since $\rho$ gives
the
weak-$*$ topology of $S$. It follows from Ascoli's theorem \cite{Ru} that
$(\Cal B_r)\ \hat{}$  is totally bounded for $\| \ \|_{\infty}$. By
Condition 1.5 this means that $\Cal B_r$ is totally
bounded for $\| \ \|$ as a subset of $A$, as desired.

For the other direction we do not need Condition 1.5. We suppose now that
the image of $\Cal L_1$ in $\tilde L$ is totally bounded for $\| \
\|^{\sptilde}$. Let $\mu\in S$ and $\varepsilon > 0$ be given, and let
$B(\mu,\varepsilon)$ be the $\rho$-ball of radius $\varepsilon$ about $\mu$
in $S$. In view of Proposition 1.4 it suffices to show that
$B(\mu,\varepsilon)$ contains a weak-$*$ neighborhood of $\mu$. Now by the
total boundedness of the image of $\Cal L_1$ we can find $a_1,\dots
,a_n\in\Cal L_1$ such that the $\| \ \|^{\sptilde}$-balls of radius
$\varepsilon/3$ about the $\hat a_j$'s cover the image of $\Cal L_1$. We
now show that the weak-$*$ neighborhood
$$
\Cal O=\Cal O(\mu,\{a_j\},\varepsilon/3)= \{\nu\in S:|(\mu-\nu)(a_j)| <
\varepsilon/3, \ 1\leq j\leq n\} $$
is contained in $B(\mu,\varepsilon)$. Consider any $a\in\Cal L_1$. There is
a $j$ and a $b\in\Cal K$, depending on $a$, such that
$$
\|a-a_j-b\| <\varepsilon/3 \ .
$$
Hence for any $\nu\in\Cal O$ we have
$$
\align
|\mu(a)-\nu(a)| & \ \leq \
|\mu(a)-\mu(a_j+b)|+|\mu(a_j+b)-\nu(a_j+b)|+ |\nu(a_j+b)-\nu(a)| \\
& \ < \ \varepsilon/3 + |\mu(a_j)-\nu(a_j)| + \varepsilon/3 \ < \
\varepsilon \ .
\endalign
$$
Thus $\rho(\mu,\nu) < \varepsilon$. Consequently $\Cal O\subseteq
B(\mu,\varepsilon)$ as desired. \qed
\enddemo

\bigskip
Examination of the proof of the above theorem suggests a reformulation
which provides a convenient subdivision of the problem of showing for
specific examples that the $\rho$-topology agrees with the weak-$*$
topology. We will use this reformulation in the next sections.

\proclaim{1.9 \ Theorem} Let the data be as in {\rm{(1.3a--e).}} Then the
$\rho$-topology on $S$ will agree with the weak-$*$ topology if the
following three hypotheses are satisfied: \roster
\item"i)" Condition 1.5 holds.
\item"ii)" $\rho$ is bounded.
\item"iii)" The set $\Cal B_1=\{a\in\Cal L: L(a)\leq 1 \ {\text{and}} \
\|a\|\leq 1\}$ is totally bounded in $A$ for $\|\ \|$. \endroster
Conversely, if Condition 1.5 holds and if the $\rho$-topology agrees with
the weak-$*$ topology, then the above three conditions are satisfied.
\endproclaim

\demo{Proof} If conditions i) and ii) are satisfied, then, just as in the
first part of the proof of Theorem 1.8, there is a constant $r$ such that
the image of $\Cal B_r$ in $\tilde\Cal L$ contains the image of $\Cal L_1$.
But $\Cal B_r\subseteq r\Cal B_1$. Thus if $\Cal B_1$ is totally bounded
then so is $\Cal B_r$, as is then the image of $\Cal L_1$. Then we can
apply Theorem 1.8 to conclude that the $\rho$-topology agrees with the
weak-$*$ topology.

Conversely, if the $\rho$-topology and the weak-$*$ topology agree, then
condition ii) holds by Proposition 1.6. But by the first part of the proof
of Theorem 1.8 there is then a constant $r$ such that $\Cal B_r$ is totally
bounded. By scaling we see that $\Cal B_1$ is also. \qed
\enddemo

\bigskip
We remark that if we take any 1-dimensional subspace $\Cal K$ of an
infinite-dimensional normed space $A$, set ${\Cal L} = A$, and let $L$ be
the pull-back to $A$ of $\|\ \|\ \tilde{}$   on $A/{\Cal K}$, we obtain an
example where $\rho$ is bounded but the image of ${\Cal L}_1$ in ${\Cal
L}\sptilde$  is not totally bounded, nor is ${\Cal B}_1$ totally bounded in
$A$.

In the next sections we will find very useful the following:

\proclaim{1.10 \ Comparison Lemma} Let the data be as in {\rm{(1.3a--e)}}.
Suppose we have a subspace $\Cal M$ of $\Cal L$ which contains $\Cal K$ and
separates the points of $S$, and a semi-norm $M$ on $\Cal M$ which takes
value $0$ exactly on $\Cal K$. Let $\rho_L$ and $\rho_M$ denote the
corresponding metrics on $S$ (possibly taking value $+\infty$). Assume that
$$
M\ \geq \ L \ \ {\text{on}} \ \ \Cal M,
$$
in the sense that $M(a)\ \geq \ L(a)$ for all $a\in\Cal M$. Then
$$
\rho_M \ \leq \ \rho_L \ ,
$$
in the sense that $\rho_M(\mu,\nu)\leq \rho_L (\mu,\nu)$ for all
$\mu,\nu\in S$. Thus
\roster
\item"i)" If $\rho_L$ is finite then so is $\rho_M$. \item"ii)" If $\rho_L$
is bounded then so is $\rho_M$. \item"iii)" If the $\rho_L$-topology on $S$
agrees with the weak-$*$ topology then so does the $\rho_M$-topology.
\endroster
\endproclaim

\demo{Proof} If $a\in\Cal M$ and $M(a)\leq 1$ then $L(a)\leq 1$. Thus the
supremum defining $\rho_M$ is taken over a smaller set than that for
$\rho_L$, and so $\rho_M\leq\rho_L$. Conclusions i) and ii) are then
obvious. Conclusion iii) follows from the fact that a continuous bijection
from a compact space to a Hausdorff space is a homeomorphism. \qed
\enddemo

\bigskip
For later use we record the following easily verified fact.

\proclaim{1.11 \ Proposition} Let data be as above. Let $t$ be a strictly
positive real number. Set $M=tL$ on $\Cal L$. Then $\rho_M=t^{-1}\rho_L$.
Thus properties for $\rho_L$ of finiteness, boundedness, and agreement of
the $\rho_L$-topology with the weak-$*$ topology carry over to $\rho_M$.
\endproclaim

\goodbreak
\bigskip\noindent{\bf 2. Metrics from actions and length functions}

\medskip\noindent
Let $G$ be a compact group (with identity element denoted by $e$). We
normalize Haar measure to give $G$ mass 1. We recall that a length function
on a group $G$ is a continuous non-negative
real-valued function, $\ell$, on $G$ such that $$
\ell(xy)\leq \ell(x)+\ell(y) \quad {\text{for}} \ \ x,y\in G, \tag 2.1a
$$
$$
\ell(x^{-1}) = \ell(x),
\tag 2.1b
$$
$$
\ell(x)=0 \quad {\text{exactly if}} \ \ x=e. \tag 2.1c
$$

Length functions arise in a number of ways. For example, if $\pi$ is a
faithful unitary representation of $G$ on a finite-dimensional Hilbert
space, then we can set $\ell(x)=\|\pi_x-\pi_e\|$. We will see another way
in the next section. We will assume for the rest of this section that a
length function has been chosen for $G$.

Let $A$ be a unital $C^*$-algebra, and let $\a$ be an action (strongly
continuous) of $G$ by automorphisms of $A$. We let $\Cal L$ denote the set
of Lipschitz elements of $A$ for $\a$ (and $\ell$), with corresponding
Lipschitz semi-norm $L$. That is \cite{Ro1, Ro2}, for $a\in A$ we set
$$
L(a)=\sup \{\|\a_x(a)-a\|/ \ell(x): x\neq e\} \ , $$
which may have value $+\infty$, and we set $$
\Cal L=\{a\in A: L(a) < \infty\} \ .
$$
It is easily verified that $\Cal L$ is a $*$-subalgebra of $A$, and that
$L$ satisfies the Leibniz property 1.2. (More generally, for $0 < r < 1$ we
could define $L^r$ by
$$ L^r(a) = \sup\{\|\a_x(a)-a\|/(\ell(x))^r:x\neq e\}$$ along the lines
considered in \cite{Ro1, Ro2}. For actions on the non-commutative torus
this has been studied in \cite{Wv2}, but we will not pursue this here.)

It is not so clear whether
$\Cal L$ is carried into itself by $\a$, but we do not need this fact here.
(For Lie groups see theorem 4.1 of \cite{Ro1} or the comments after theorem
6.1 of \cite{Ro2}.) Let us consider, however, the $\a$-invariance of $L$.
We find that

$$
\align
L(\a_z(a))&= \sup \{\|\a_z(\a_{z^{-1}xz}(a)-a)\|/ \ell(x): x\neq e\} \\
&= \sup
\{\|\a_x(a)-a\|/ \ell(zxz^{-1}): x\neq e\}.
\endalign
$$
Thus if $\ell(zxz^{-1}) = \ell(x)$ for all $x, z \in G$, then $L$ is
$\a$-invariant, and $\Cal L$ is carried into itself by $\a$. The metric
$\rho$ on $S$ defined by $L$ will then be $\a$-invariant for the evident
action on $S$. But we will not discuss this matter further here.

\proclaim{2.2 \ Proposition} The $*$-algebra $\Cal L$ is dense in $A$.
\endproclaim

\demo{Proof} For $f \in L^1(G)$ we define $\a_f$ as usual by $\a_f(a) =
\int f(x)\a_x(a)\ dx$. It is standard \cite{BR} that as $f$ runs through an
``approximate delta-function'', $\a_f(a)$ converges to $a$. Thus the set of
elements of form $\a_f(a)$ is dense in $A$. Let $\lambda$ denote the action
of $G$ by left translation of functions on $G$. A quick standard
calculation shows that $\a_x(\a_f(a)) = \a_{\lambda_x(f)}(a)$. Thus
$$
\|\a_x(\a_f(a)) - \a_f(a)\| = \|\a_{(\lambda_xf -f)}(a)\| \leq \|\lambda_xf
-f\|_1 \|a\|,
$$
where $\|\ \|_1$ denotes the usual $L^1$-norm. Thus we see that $\a_f(a)
\in \Cal L$ if $f \in Lip^1_\lambda$, the space of Lipschitz functions in
$L^1(G)$ for $\lambda$ (and $\ell$).

Consequently it suffices to show that $Lip^1_\lambda$ is dense in $L^1(G)$.
We first note that it contains a non-trivial element, namely $\ell$ itself.
For if $x, y \in G$, then $$
|(\lambda_x\ell)(y) - \ell(y)| = |\ell(x^{-1}y) - \ell(y)| \leq \ell(x),
$$
where the inequality follows from 2.1a and 2.1b above. We momentarily
switch attention to $C(G)$ with $\|\ \|_\infty$, and the action $\lambda$
of $G$ on it. Of course $\ell \in C(G)$. The above inequality then says
that $\ell \in Lip^\infty_\lambda$, the space of Lipschitz functions in
$C(G)$ for $\lambda$. But as mentioned earlier, $Lip^\infty_\lambda$ is
easily seen to be a $*$-subalgebra of $C(G)$ for the pointwise product, and
it contains the constant functions. Furthermore, a simple calculation shows
that $Lip^\infty_\lambda$ is carried into itself by {\it right}
translation. Since $Lip^\infty_\lambda$ contains $\ell$, which separates
$e$ from any other point, it follows that $Lip^\infty_\lambda$ separates
the points of $G$. Thus $Lip^\infty_\lambda$ is dense in $C(G)$ by the
Stone-Weierstrass theorem. Since $\|\ \|_\infty$ dominates $\|\ \|_1$ for
compact $G$, it follows that $Lip^1_\lambda$ is dense in $L^1(G)$ as
needed. \qed
\enddemo

For simplicity of exposition we will deal only with the case in which we
obtain metrics on the entire state space of the $C^*$-algebra A. For this
purpose we want the subspace where $L$ takes the value 0 to be
one-dimensional. It is evident that $L$ takes value 0 on exactly those
elements of $A$ which
are $\a$-invariant, and in particular on the scalar multiples of the
identity element of $A$. Thus we need to assume that the action $\a$ is
{\it ergodic}, in the sense that the only $\a$-invariant elements are the
scalar multiples of the identity.

The main
theorem of this section is:

\proclaim{2.3 \ Theorem} Let $\a$ be an ergodic action of a compact group
$G$ on a unital $C^*$-algebra $A$. Let $\ell$ be a length function on $G$,
and define $\Cal L$ and $L$ as above. Let $\rho$ be the corresponding
metric on the state space $S$ of $A$. Then the $\rho$-topology on $S$
agrees with the weak-$*$ topology.
\endproclaim

\demo{Proof} Because $\Cal L$ is dense by Proposition 2.2, it separates the
points of $S$. Consequently the conditions 1.3a--e are fulfilled (for the
evident $\eta$). Thus $L$ indeed defines a metric, $\rho$, on $S$ (perhaps
taking value $+\infty$).

Since $G$ is compact, we can average $\a$ over $G$ to obtain a conditional
expectation from $A$ onto its fixed-point subalgebra. Because we assume
that $\a$ is ergodic, this conditional expectation can be viewed as a state
on $A$. By abuse of notation we will denote it again by $\eta$, since it
extends the evident state $\eta$ on the fixed-point algebra. Thus
$$
\eta(a)=\int_G \a_x(a) \ dx
$$
for $a\in A$, interpreted as a complex number when convenient.

We will follow the approach suggested by Theorem 1.9. Now hypothesis (i) of
that theorem is
satisfied in the present setting, as discussed right after Condition 1.5
above. We now check hypothesis (ii), that is: \enddemo

\proclaim{2.4 \ Lemma} $\rho$ is bounded. \endproclaim

\demo{Proof} Let $\mu\in S$. Then for any $a\in\Cal L$ we have $$
|\mu(a)-\eta(a)| = |\int\mu(a)dx-\mu(\int \a_x(a)dx)| =
|\int\mu(a-\a_x(a))dx| \leq L(a)\int_G \ell(x)dx \ . $$
It follows that $\rho(\mu,\eta)\leq\int\ell(x)dx$. Thus for any $\mu,\nu\in
S$ we have
$$
\rho(\mu,\nu)\leq 2\int_G \ell(x)dx \ ,
$$
which is finite since $\ell$ is bounded. \qed \enddemo

\bigskip
We now begin the verification of hypothesis (iii) of Theorem 1.9. For this
we need the unobvious fact \cite{HLS, Bo} that because $G$ is compact and
$\a$ is ergodic, each irreducible representation of $G$ occurs with at most
finite multiplicity in $A$. (In \cite{HLS} it is also shown that $\eta$ is
a trace, but we do not need this fact here.) The following lemma is
undoubtedly well-known, but I do not know a reference for it.

\proclaim{2.5 \ Lemma} Let $\a$ be a (strongly continuous) action of a
compact group $G$ on a Banach space $A$. Suppose that each irreducible
representation of $G$ occurs in $A$ with at most finite multiplicity. Then
for any $f\in L^1(G)$ the operator $\a_f$ defined by
$$
\a_f(a)=\int_G f(x)\a_x(a)dx
$$
is compact.
\endproclaim

\demo{Proof} If $f$ is a coordinate function for an irreducible
representation $\pi$ of $G$, then it is not hard to see (ch. IX of
\cite{FD}) that $\a_f$ will have range in the $\pi$-isotypic component of
$A$, which we are assuming is finite-dimensional. Thus $\a_f$ is of finite
rank in this case. But by the Peter-Weyl theorem \cite{FD} the linear span
of the coordinate functions for all irreducible representations is dense in
$L^1(G)$. So any $\a_f$ can be approximated by finite rank operators. \qed
\enddemo

\bigskip
\demo{Proof of Theorem 2.3} We show now that $\Cal B_1$, as in (iii) of
Theorem 1.9, is totally bounded. Let $\varepsilon> 0$ be given. Since
$\ell(e)=0$ and $\ell$ is continuous at $e$, we can find $f\in L^1(G)$ such
that $f\geq 0$,
\ $\int_G f(x)dx=1$, and $\int_G f(x)\ell(x)dx < \varepsilon/2$. By the
previous lemma $\a_f$ is compact. Since $\Cal B_1$ is bounded, it follows
that $\a_f(\Cal B_1)$ is totally bounded. Thus it can be covered by a
finite number of balls of radius $\varepsilon/2$. But for any $a\in\Cal
B_1$ we have $$
\align
\|a-\a_f(a)\| & = \|a\int f(x)dx-\int f(x)\a_x(a)dx\| \leq \int
f(x)\|a-\a_x(a)\|dx \\ &\leq
L(a)\int f(x)\ell(x)dx
\leq \varepsilon/2 \ .
\endalign
$$
Thus $\Cal B_1$ itself can be covered by a finite number of balls of radius
$\varepsilon$. \qed
\enddemo

\goodbreak
\bigskip\bigskip\noindent
{\bf 3. Metrics from actions of Lie groups}

\medskip\noindent
We suppose now that $G$ is a connected Lie group (compact). We let $\g$
denote the Lie algebra of $G$. Fix a norm $\|\ \|$ on $\g$. For any action
$\a$ of $G$ on a Banach space $A$ we let $A^1$ denote the space of
$\a$-differentiable elements of $A$. Thus \cite{BR} if $a\in A^1$ then for
each $X\in\g$ there is a $d_Xa\in A$ such that
$$
\lim_{t\to 0} (\a_{\exp(tX)}(a)-a)/t=d_Xa \ , $$
and $X\mapsto d_Xa$ is a linear map from $\g$ into $A$, which we denote by
$da$. Since $\g$ and $A$ both have norms, the operator norm, $\|da\|$, of
$da$ is defined (and finite). A standard smoothing argument \cite{BR}
shows that $A^1$ is dense in $A$.

Suppose now that $A$ is a $C^*$-algebra and that $\a$ is an action by
automorphisms of $A$.
We can set $\Cal L=A^1$ and $L(a)=\|da\|$. It is easily verified that $\Cal
L$ is a $*$-subalgebra of $A$ and that $L$ satisfies the Leibniz property
1.2, though we do not need these facts here. Because $G$ is connected,
$L(a)=0$ exactly if $a$ is $\a$-invariant.

\proclaim{3.1 \ Theorem} Let $G$ be a compact connected Lie group, and fix
a norm on $\g$.
Let $\a$ be an ergodic action of $G$ on a unital $C^*$-algebra $A$. Let
$\Cal L=A^1$ and $L(a)=\|da\|$, and let $\rho$ denote the corresponding
metric on the state space $S$. Then the $\rho$-topology on $S$ agrees with
the weak-$*$ topology. \endproclaim

\demo{Proof} Choose an inner-product on $\g$. Its corresponding norm is
equivalent to the given norm, and so by the Comparison Lemma 1.10 it
suffices to deal with the norm from the inner-p
roduct. We can
left-translate this inner-product over $G$ to obtain a left-invariant
Riemannian metric on $G$, and then a corresponding left-invariant ordinary
metric on $G$. We let $\ell(x)$ denote the corresponding distance from
$x$ to $e$. Then $\ell$ is a continuous length function on $G$ satisfying
conditions 2.1 \cite{G, Ro2}.

Then the elements of $\Cal L = A^1$ are Lipschitz for $\ell$. This
essentially just involves the following standard argument \cite{G, Ro2},
which we include for the reader's convenience. Let $a\in A^1$ and let $c$
be a smooth path in $G$ from $e$ to a point $x\in G$. Then $\phi$, defined
by $\phi(t)=\a_{c(t)}(a)$, is differentiable, and so we have $$
\|\a_x(a)-a\| =\|\int\phi'(t)dt\| \leq
\int \|\a_{c(t)} (d_{c'(t)}a)\|dt \leq \|da\| \int \|c'(t)\|dt \ .
$$
But the last integral is just the length of $c$. Thus from the definition
of the ordinary metric on $G$, with its length function $\ell$, we obtain
$$
\|\a_x(a)-a\|\leq \|da\|\ell(x) \ .
$$
(Actually, the above argument works for any norm on $\g$.) Then if we let
$\Cal L_0$ and $L_0$ be defined just in terms of $\ell$ as in the previous
section, we see that $\Cal L\subseteq \Cal L_0$ and $L_0\leq L$. Thus we
are exactly in position to apply the Comparison Lemma 1.10 to obtain the
desired conclusion. \qed
\enddemo

\bigskip
We remark that Weaver (theorem 24 of \cite{Wv1}) in effect proved for this
setting
the total boundedness of $\Cal B_1$ for the particular case of
non-commutative 2-tori, by different methods.

\bigskip\bigskip\noindent{\bf 4. Metrics from Dirac operators}

\medskip\noindent
Suppose again that $G$ is a compact connected Lie group, and that $\a$ is
an ergodic action of $G$ on a unital $C^*$-algebra $A$. Let $\g$ denote the
Lie algebra of $G$, and let $\g'$ denote its vector-space dual. Fix any
inner-product on $\g'$. We will denote it by $g$, or by $\< \ , \>_g$, to
distinguish it from the Hilbert space inner-products which will arise.

With this data we can define a spectral triple \cite{C1, C2, C3} for
$A$. For simplicity of exposition we will not include gradings and real
structure, and we will oversimplify our treatment of spinors, since the
details are not essential for our purposes. But with more care they can be
included. (See, e.g. \cite{V, VB}.) We proceed
as follows. Let $C={\op{Clif}}(\g',-g)$ be the complex Clifford
$C^*$-algebra over $\g'$ for $-g$. Thus each $\omega\in\g'$ determines a
skew-adjoint element of $C$ such that $$
\omega^2=-\<\omega,\omega\>_g \ 1_C \ .
$$
Depending on whether $\g$ is even or odd dimensional, $C$ will be a full
matrix algebra, or the direct sum of two such. We let $\Cal S$ be the
Hilbert space of a finite-dimensional faithful representation of $C$ (the
``spinors").

Let $A^{\infty}$ denote the space of smooth elements of $A$. (We could just
as well use the $A^1$ of the previous section. We use $A^{\infty}$ here for
variety. It is still a dense $*$-subalgebra \cite{BR}.) Let
$W=A^{\infty}\otimes\Cal S$, viewed as a free right $A^{\infty}$-module.
From the Hilbert-space inner-product on $\Cal S$ we obtain
an $A^{\infty}$-valued
inner-product on $W$. Let $\eta$ be as in the previous section,
viewed as a faithful state on $A$. Combined with the $A$-valued
inner product on $W$, it gives an ordinary inner-product on $W$. We
will denote the Hilbert space completion by $L^2(W,\eta)$.

Now $A^{\infty}$ and $C$ have evident commuting left actions on $W$. These
are easily seen to give $*$-representations of $A$ and $C$ on
$L^2(W,\eta)$, which we denote by $\lambda$ and $c$ respectively.

We define the Dirac operator, $D$, on $L^2(W,\eta)$ in the usual way. Its
domain will be $W$, and it is defined as the composition of operators
$$
W \overset d\to\longrightarrow
\g'\otimes W \overset i\to\longrightarrow C\otimes W \overset
c\to\longrightarrow W \ . $$
Here $d$ is the operator which takes $b\in A^{\infty}$ to $db\in\g'\otimes
A^{\infty}$, defined by $db(X)=d_X(b)$, which we then extend to $W$ so that
it takes $b\otimes s$ to $db\otimes s$. The operator $i$ just comes from
the canonical inclusion of $\g'$ into $C$. The operator $c$ just comes from
applying the representation of $C$ on $\Cal S$, and so on $W$.

It is easily seen that $D$ is a symmetric operator on $L^2(W,\eta)$. It
will not be important for us to verify that $D$ is essentially
self-adjoint, and that its closure has compact resolvant.

Let $\{e_j\}$ denote an orthonormal basis for $\g'$, and let $\{E_j\}$
denote the dual basis for $\g$. Then in terms of these bases we have
$$
D(b\otimes s)=\sum \a_{E_j} (b)\otimes c(e_j)s \ . $$
When we use this to compute $[D,\lambda_a]$ for $a\in A^{\infty}$, a
straightforward calculation shows that we obtain
$$
[D,\lambda_a](b\otimes s) = \sum (\a_{E_j} (a)\otimes c(e_j))) (b\otimes s) \ .
$$
That is,
$$
[D,\lambda_a] = \sum \a_{E_j} (a)\otimes e_j \ , \tag 4.1 $$
acting on $L^2(W,\eta)$ through the representations $\lambda$ and $c$. It
is clear from (4.1) that $[D,\lambda_a]$ is bounded for the operator norm
from $L^2(W,\eta)$.

We can now set $\Cal L=A^{\infty}$, and
$$
L(a)=\|[D,\lambda_a]\| \ .
$$
It is clear that $L(1_A)=0$. To proceed further we compare $L$ with the
semi-norm of the last section. If we view $\g'$ as contained in the
$C^*$-algebra $C$, we have $e^2_j=-1$ and $e^*_j=-e_j$ for each $j$. In
particular, $\|e_j\|=1$. From (4.1) it is then easy to see that there is a
constant, $K$, such that $$
L(a)\leq K\|da\|
$$
for all $a\in\Cal L$, where $\|da\|$ is as in the previous section, for the
inner-product dual to that on $\g'$. However, what we need is an inequality
in the reverse direction so that we will be able to apply the Comparison
Lemma 1.10.

For this purpose, consider any element $t=\sum b_j\otimes e_j$ in $A\otimes
C$, with the $e_j$ as above. Let $f_j=ie_j$, so that $f^*_j=f_j$, \
$f^2_j=1$, and $f_jf_k=-f_kf_j$ for $j\neq k$. Let $p_j=(1+f_j)/2$ and
$q_j=1-p_j=(1-f_j)/2$, both being self-adjoint projections. Then
$p_jf_k=f_kq_j$ for $j\neq k$. Consequently $p_jf_kp_j=0=q_jf_kq_j$ for
$j\neq k$. Thus $$
(1\otimes p_j)t(1\otimes p_j) = b_j\otimes p_je_jp_j = b_j\otimes ip_j
$$
and
$$
(1\otimes q_j)t(1\otimes q_j) = -b_j\otimes iq_j \ . $$
Since at least one of $p_j$ and $q_j$ must be non-zero, it becomes clear that
$\|t\| \geq \|b_j\|$
for each $j$. When we apply this to (4.1) we see that $$
L(a)\geq \|\a_{E_j}(a)\|
$$
for each $j$. Consequently, for a suitable constant $k$ we have $$
L(a)\geq k\|da\| \ ,
$$
where again $\|da\|$ is as in the previous section. On applying Proposition
1.11,
Theorem 3.1, and the Comparison Lemma 1.10, we obtain the proof of:

\proclaim{4.2 \ Theorem} Let $\a$ be an ergodic action of the compact
connected Lie group $G$ with Lie algebra $\g$ on the unital $C^*$-algebra
$A$. Pick any inner-product on the dual, $\g'$, of $\g$. Let $D$ denote the
corresponding Dirac operator, as defined above. Let $\Cal L=A^{\infty}$,
and let $L$ be defined by $$
L(a)=\|[D,a]\|
$$
for $a\in A$. Let $\rho$ be the corresponding metric on $S$. Then the
$\rho$-topology on $S$ agrees with the weak-$*$ topology. \endproclaim

\enddocument